\documentclass[10pt,twoside,leqno]{amsart}

\def\C{\mathbb C}
\def\R{\mathbb R}
\def\N{\mathbb N}

\def\cK{{\mathcal K}}
\def\cH{{\mathcal H}}
\def\fH{{\mathfrak H}}

\def\Mt{{\widetilde{M}}}
\def\mt{{\widetilde{m}}}
\def\Bt{{\widetilde{B}}}
\def\At{{\widetilde{A}}}
\def\Gammat{{\widetilde{\Gamma}}}
\def\gammat{{\widetilde{\gamma}}}

\def\cL{{\mathcal L}}

\def\eps{\varepsilon}
\def\wlim{w-\lim}

\def\ker{{\mathrm{ker\,}}}
\def\Ran{{\mathrm{Ran\,}}}
\def\rg{{\mathrm{Ran\,}}}
\def\dive{{\mathrm{div\,}}}

\def\Im{{\mathrm{Im}}}
\def\into{{\int_\Omega}}
\def\intdo{{\int_{\partial\Omega}}}

\newcommand{\ie}{i.e.~}
\newcommand{\cf}{c.f.~}

\newcommand{\norm}[1]{\left\Vert#1\right\Vert}
\newcommand{\trinorm}[1]{\left\vert\!\vert\! \vert#1\vert\!\vert\!\right\vert}

\newcommand{\be}{\begin{equation}}
\newcommand{\ee}{\end{equation}}

\newcommand{\beq}{\begin{equation}}
\newcommand{\enq}{\end{equation}}

\def\llangle{\left\langle}
\def\rrangle{\right\rangle}
\def\({\left(}
\def\){\right)}

\usepackage{latexsym}

\usepackage{amsmath}
\usepackage{amssymb}
\usepackage{amsthm}
\usepackage{latexsym}
\usepackage[centertags]{amsmath}

\theoremstyle{definition}
\theoremstyle{plain}

\newtheorem{proposition}{Proposition}[section]
\newtheorem{theorem}[proposition]{Theorem}
\newtheorem{lemma}[proposition]{Lemma}
\newtheorem{corollary}[proposition]{Corollary}
\newtheorem{remark}[proposition]{Remark}
\newtheorem{definition}[proposition]{Definition}
\newtheorem{example}[proposition]{Example}

\numberwithin{equation}{section}

\allowdisplaybreaks

\begin{document}

\title[Boundary triplets and $M$-functions]{Boundary triplets and $M$-functions for non-selfadjoint
 operators, with applications to elliptic PDEs and block operator
 matrices}

\author{Malcolm Brown}
\address{School of Computer Science, Cardiff University, Queen's Buildings, 5 The Parade, Cardiff CF24 3AA, UK}
\email{Malcolm.Brown@cs.cardiff.ac.uk}

\author{Marco Marletta}
\address{School of Mathematics, Cardiff University, Senghennydd Road, Cardiff CF24 4AG, UK}
\email{MarlettaM@cardiff.ac.uk}

\author{Serguei Naboko}
\address{Department of Math.~Physics, Institute of Physics, St.~Petersburg State University, 1 Ulianovskaia, St.~Petergoff, St.~Petersburg, 198504, Russia}
\email{naboko@snoopy.phys.spbu.ru}

\author{Ian Wood}
\thanks{Serguei Naboko wishes to thank British EPSRC for supporting his visit to Cardiff under the grant EP/C008324/1 "Spectral Problems on Families of Domains and Operator M-functions". He also wishes to thank C
ardiff University for hospitality during the visit. Ian Wood wishes to thank British EPSRC for support under the same grant.\\
The authors are grateful to Professors Yury Arlinskii, Fritz Gesztesy, Gerd Grubb, Mark Malamud, Andrea Posilicano and Vladimir Ryzhov for useful comments on an earlier version of this paper.}
\address{Institute of Mathematical and Physical Sciences, University of Wales Aberystwyth, Penglais, Aberystwyth, Ceredigion SY 23 3BZ, UK}
\email{ian.wood@aber.ac.uk}


\begin{abstract}
Starting with an adjoint pair of operators, under suitable abstract versions of standard PDE hypotheses, we consider the Weyl $M$-function of extensions of the operators. The extensions are  determined by abstract boundary conditions  and we establish results on the relationship between the $M$-function as an analytic function of a spectral parameter and the spectrum of the extension. We also give an example where the $M$-function does not contain the whole spectral information of the resolvent, and show that the results can be applied to elliptic PDEs where the $M$-function corresponds to the Dirichlet to Neumann map.
\end{abstract}

\maketitle

\section{Introduction}
The theory of boundary value spaces associated with symmetric operators
has its origins in the work of Ko\v{c}ube\u{\i} \cite{Ko75} and Gorbachuk and Gorbachuk
\cite{Gorbachuk} and has been the subject of intense activity in the former 
Soviet Union, with major contributions from many authors. While we cannot
undertake a comprehensive survey of the literature here, we recommend that the
reader consult the works of Derkach and Malamud who developed the theory of the Weyl-$M$-function in the context of boundary value spaces (e.g. \cite{DM87,Derkach}); the work of V.A. Mikhailets (e.g. the very elegant application of the theory of boundary value 
spaces by Mikhailets and Sobolev \cite{MikhailetsS} to the common eigenvalue 
problem for periodic Schr\"{o}dinger operators); the work of Kuzhel and Kuzhel 
(e.g. \cite{Kuzhel1,Kuzhel2}); the work of Brasche, Malamud and Neidhardt 
(e.g. \cite{Brasche}); the work of Storozh (in particular, \cite{Storozh}) and the recent work of Kopachevski\u\i\ and Kre\u\i n \cite{KK04} and Ryzhov \cite{Ryz07} on abstract Green's formulae, again Ryzhov \cite{Ryzhov} on functional models and Posilicano \cite{Pos07} characterising extensions and giving some applications to PDEs.

Adjoint pairs of second order elliptic operators, their extensions and boundary value problems were studied in the paper of Vishik \cite{Vis52}. For adjoint pairs of abstract operators, boundary triplets were introduced by Vainerman \cite{Vai80} 
and Lyantze and Storozh \cite{LS83}. Many of the results proved for the symmetric case, such 
as characterising extensions of the operators and investigating spectral properties via the 
Weyl-$M$-function, have subsequently been extended for this situation: see, for instance,
Malamud and Mogilevski \cite{MMM1} for adjoint pairs of operators, Langer and Textorius
\cite{LT} and Malamud \cite{M} for adjoint pairs of contractions, and Malamud and Mogilevski
\cite{MMM2,MM02} for adjoint pairs of linear relations. 
For the case of sectorial operators and their $M$-functions we should mention
especially the work of Arlinskii \cite{Arl99,Arl00,Arl01} 
who uses sesquilinear form methods. 
The approach using adjoint pairs of operators does not
require any assumption that the operators be sectorial. The price which must be paid
for this is that there are other hypotheses (e.g. non-emptiness of the resolvent set of 
certain operators or, in our approach, an abstract unique continuation assumption)
which must be verified before this approach can be applied.

In the context of PDEs there has also been extensive work on Dirichlet to Neumann maps, also 
sometimes known as Poincar\'{e}-Steklov operators, especially in the inverse problems 
literature. These operators have physical meaning, associating, for instance, a surface 
current to an applied voltage. For some applications of them to quantum networks we refer 
to recent papers by Pavlov et al.~\cite{HarmerPavlovY} and \cite{Pavlov}. These maps 
are, in some sense, the natural PDE realization of the abstract $M$-function 
which appears in the theory of boundary value spaces. Amrein and Pearson \cite{AP04} generalised several results from the classical Weyl-$m$-function for the one-dimensional Sturm-Liouville problem to the case of Schr\"odinger operators, calling them $M$-functions, in particular they were able to show nesting results for families of $M$-functions on exterior domains.
However there have been relatively few applications of the theory of boundary value spaces to PDEs. A chapter in Gorbachuck and Gorbachuk \cite{Gorbachuk} deals with
a PDE on a tubular domain by reduction to a system of ODEs with operator
coefficients, and there are some papers which deal with special perturbations of
PDE problems which result in symmetric operators with (crucially) finite
deficiency indices, e.g. the very recent paper of Br\"{u}ning,
Geyler and Pankrashkin \cite{Geyler}. The case of symmetric operators with infinite deficiency indices is studied by Behrndt and Langer in \cite{BL07}. However for symmetric elliptic PDEs a 
concrete realization of the boundary value operators whose existence is 
guaranteed by the abstract theory, and a precise description of the relationship 
between the abstract $M$-function and the classical Dirichlet to Neumann map,  
requires a technique due to Vishik \cite{Vis52} and Grubb \cite{Grubb} in the choice of the 
boundary value operators 
which we describe in this paper.

In this paper we consider the non-symmetric case. Using the setting of boundary triplets from Lyantze and Storozh
\cite{LS83}, we introduce an $M$-function and prove
the following results:
\begin{itemize}
\item[i.] the relationship between poles of the $M$-function as an analytic 
function of a spectral parameter and eigenvalues of a corresponding operator 
determined by abstract boundary conditions, under a new abstract unique
continuation hypothesis which is natural in the context of PDEs;
\item[ii.] results concerning behaviour of the $M$-function near the 
 essential spectrum;
\item[iii.] a proof that the $M$-function does not contain the whole spectral information of the resolvent, by consideration of a 
 Hain-L\"{u}st problem;
\item[iv.] results concerning the analytic behaviour of Dirichlet to Neumann
 maps for elliptic PDEs, though these have also been obtained recently in a 
 concrete way by F.~Gesztesy et al.~\cite{GMZ07}.
\end{itemize}

\section{Basic concepts and notation}

Throughout, we will make the following assumptions:
\begin{enumerate}
  \item $A$ and $\At$ are closed densely defined operators on a Hilbert space $H$.
  \item $A$ and $\At$ are an adjoint pair, \ie $A^*\supseteq\At$ and $\At^*\supseteq A$.
  \item Whenever considering $D(\At^*)$ as a linear space it will be equipped with the graph norm. Since $\At^*$ is closed, this
makes $D(\At^*)$ a Hilbert space.
\end{enumerate}

\begin{proposition}(Lyantze, Storozh '83). For each adjoint pair of closed densely defined operators on $H$,
there exist ``boundary spaces'' $\cH$, $\cK$ and ``boundary operators''
  \[ \Gamma_1:D(\At^*)\to\cH,\quad \Gamma_2:D(\At^*)\to\cK,\quad \Gammat_1:D(A^*)\to\cK\quad \hbox{ and }\quad  \Gammat_2:D(A^*)\to\cH \]
  such that for $u\in D(\At^*) $ and $v\in D(A^*)$ we have an abstract Green formula
        \begin{equation}\label{Green}
          (\At^* u, v)_H - (u,A^*v)_H = (\Gamma_1 u, \Gammat_2 v)_\cH - (\Gamma_2 u, \Gammat_1v)_\cK.
        \end{equation} 
  The boundary operators $\Gamma_1$, $\Gamma_2$, $\Gammat_1$ and $  \Gammat_2 $ are bounded with respect to the graph norm and surjective. Moreover, we have 
\begin{equation}\label{domains}
	D(A)= D(\At^*)\cap\ker\Gamma_1\cap \ker\Gamma_2 \quad \hbox{ and } \quad D(\At)= D(A^*)\cap\ker\Gammat_1\cap \ker\Gammat_2.
\end{equation}
The collection $\{\cH\oplus\cK, (\Gamma_1,\Gamma_2), (\Gammat_1,\Gammat_2)\}$ is called a boundary triplet for the adjoint pair $A,\At$.
\end{proposition}

\begin{proof}
The proof in Russian is in \cite[Chapter 4]{LS83}. For the more general situation of linear relations a proof in English can be found in \cite[Section 3.2]{MM02}. 
\end{proof}

\begin{remark}
Using this setting, in \cite{MM02} Malamud and Mogilevskii go on to define Weyl $M$-functions and $\gamma$-fields associated with boundary triplets and to obtain Kre\u{\i}n formulae for the resolvents. In the same spirit we introduce $M$-functions and what we call the solution operator. In our setting, these will depend on a parameter given by an operator $B\in\cL(\cK,\cH)$.  To take account of this technical difference and to keep this paper as self-contained as possible we will develop the full theory in Sections 2 and 3 here, noting that similar definitions and results can be found in \cite{MM02}. 
\end{remark}

\begin{definition}\label{defmfn}
We consider the following extensions of $A$ and $\At$: Let $B\in\cL(\cK,\cH)$ and $\Bt\in\cL(\cH,\cK)$ and define
\[ A_B:=\At^*\vert_{\ker(\Gamma_1-B\Gamma_2)} \hbox{ and } \At_\Bt:=A^*\vert_{\ker(\Gammat_1-\Bt\Gammat_2)}.\]  In the following, we will always assume $\rho(A_B)\neq\emptyset$, in particular $A_B$ will be a closed operator.

For $\lambda\in\rho(A_B)$, we define the $M$-function via
\[ M_B(\lambda):\Ran(\Gamma_1-B\Gamma_2)\to\cK,\ M_B(\lambda)(\Gamma_1-B\Gamma_2) u=\Gamma_2 u \hbox{ for all } u\in \ker(\At^*-\lambda)\]
and for $\lambda\in\rho(\At_\Bt)$, we define 
\[ \Mt_\Bt(\lambda):\Ran(\Gammat_1-\Bt\Gammat_2)\to\cH,\ \Mt_\Bt(\lambda)(\Gammat_1-\Bt\Gammat_2) v=\Gammat_2 v 
\hbox{ for all } v\in \ker(A^*-\lambda).\]
\end{definition}

\begin{lemma}
$M_B(\lambda)$ and $\Mt_\Bt(\lambda)$ are well-defined. 
\end{lemma}

\begin{proof}
We prove the statement for $M_B(\lambda)$. Suppose $f\in\Ran(\Gamma_1-B\Gamma_2)$, then 
there exists $u\in\ker(\At^*-\lambda)$ such that $(\Gamma_1-B\Gamma_2)u=f$. To see this, choose any $w\in D(\At^*)$ such that $(\Gamma_1-B\Gamma_2)w=f$. Let $v=-(A_B-\lambda)^{-1}(\At^*-\lambda)w$. Then $u=v+w\in\ker(\At^*-\lambda)$ and $(\Gamma_1-B\Gamma_2)(v+w)=(\Gamma_1-B\Gamma_2)w=f$.
Now assume $(\Gamma_1-B\Gamma_2) u=(\Gamma_1-B\Gamma_2) v$ for some $u,v\in\ker(\At^*-\lambda)$. 
Then $u-v\in \ker(\At^*-\lambda)\cap D(A_B)$. As $\lambda\in\rho(A_B)$, there exists $w\in H$ such that $u-v=(A_B-\lambda)^{-1}w$.
Then $0=(\At^*-\lambda)(u-v)=(\At^*-\lambda)(A_B-\lambda)^{-1}w=w$, so $u=v$, in particular, $\Gamma_2 u=\Gamma_2 v$.
\end{proof}

\section{ The solution operator $S_{\lambda,B}$}

\begin{definition}
For $\lambda\in\rho(A_B)$, we define the operator $S_{\lambda,B}:\Ran(\Gamma_1-B\Gamma_2)\to \ker(\At^*-\lambda)$ by
\begin{eqnarray}\label{slamdef}
  (\At^*-\lambda)S_{\lambda,B} f=0,\ (\Gamma_1-B\Gamma_2)S_{\lambda,B} f=f,
\end{eqnarray}
\ie $S_{\lambda,B}=\left( (\Gamma_1-B\Gamma_2)\vert_{\ker(\At^*-\lambda)}\right)^{-1}$.
\end{definition}

\begin{lemma}\label{slamwd}
$S_{\lambda,B}$ is well-defined for $\lambda\in\rho(A_B)$.
\end{lemma}

\begin{proof}
For $f\in \Ran(\Gamma_1-B\Gamma_2)$, choose any $w\in D(\At^*)$ such that $(\Gamma_1-B\Gamma_2)w=f$. Let 
$v=-(A_B-\lambda)^{-1}(\At^*-\lambda)w$. Then $v+w\in\ker(\At^*-\lambda)$ and $(\Gamma_1-B\Gamma_2)(v+w)=(\Gamma_1-B\Gamma_2)w=f$, 
so a solution to (\ref{slamdef}) exists and is given by 
\begin{equation*}
  S_{\lambda,B} f= \left( I-(A_B-\lambda)^{-1}(\At^*-\lambda)\right)w
\end{equation*}
for any $w\in D(\At^*)$ such that $(\Gamma_1-B\Gamma_2)w=f$.

Moreover, the solution to (\ref{slamdef}) is unique: Suppose $u_1$ and $u_2$ are two solutions. Then 
$(u_1-u_2)\in\ker(\At^*-\lambda)\cap \ker(\Gamma_1-B\Gamma_2)$, so $u_1-u_2\in D(A_B)$ and $(A_B-\lambda)(u_1-u_2)=0$.
 As $\lambda\in\rho(A_B)$, $u_1=u_2$. 
\end{proof}

\begin{proposition}\label{slamanaprop}
Let $f\in\Ran(\Gamma_1-B\Gamma_2)$. The map from $\rho(A_B)\to H$ given by $\lambda\mapsto S_{\lambda,B} f$ is analytic.
\end{proposition}

\begin{proof}
  Fix $\lambda_0\in\rho(A_B)$. Now choose $w=S_{\lambda_0,B}f$ in the proof of Lemma \ref{slamwd}. Then
\begin{equation}\label{slamanal}
  S_{\lambda,B} f= \left( S_{\lambda_0,B}-(A_B-\lambda)^{-1}(\At^*-\lambda)S_{\lambda_0,B}\right) f 
             = S_{\lambda_0,B}f + (\lambda-\lambda_0)(A_B-\lambda)^{-1}S_{\lambda_0,B}f,
\end{equation}
which is analytic in $\lambda$.
\end{proof}

\begin{lemma}\label{slamclosed}
Let $F:=\ker(\At^*-\lambda)$, $E:=\Ran (\Gamma_1-B\Gamma_2)$ and 
\[ \norm{u}_F^2:=\norm{u}_H^2+\norm{\At^* u}_H^2, \quad  \norm{f}_E^2:=\norm{f}_\cH^2+\norm{S_{\lambda,B} f}_F^2.\]
Then $E$ and $F$ are Hilbert spaces and the operator $S_{\lambda,B}$ with $D(S_{\lambda,B})=E\subseteq\cH$ is closed as an operator from 
$\cH$ to $D(\At^*)$.
\end{lemma}

\begin{proof}
Obviously, $\norm{\cdot}_E$ and  $\norm{\cdot}_F$ are norms induced by scalar products. It remains to prove completeness.
Since $(\At^*-\lambda):D(\At^*)\to H$ is continuous, $F$ is a closed subspace of $D(\At^*)$, hence complete.

Assume $(f_n)_{n\in\N}$ is a Cauchy sequence in $E$. Then $(f_n)_{n\in\N}$ is Cauchy in $\cH$ and converges to $f\in\cH$ and 
$(S_{\lambda,B} f_n)_{n\in\N}$ is Cauchy in $F$ and converges to $u\in F$. As $\Gamma_1-B\Gamma_2$ is continuous in the graph norm and 
$S_{\lambda,B}^{-1}:F\to E$ is given by $\Gamma_1-B\Gamma_2$, we have
\begin{eqnarray*}
  \norm{(\Gamma_1-B\Gamma_2)u -f}_\cH &=& \norm{(\Gamma_1-B\Gamma_2)u -S_{\lambda,B}^{-1}S_{\lambda,B} f_n +f_n -f}_\cH\\
  &\leq&  \norm{\Gamma_1-B\Gamma_2}_{F\to \cH} \norm{u -S_{\lambda,B} f_n}_F +\norm{f_n -f}_\cH \to 0,
\end{eqnarray*}
so $(\Gamma_1-B\Gamma_2)u =f$, \ie $f\in E$ and $S_{\lambda,B} f=u$.

Therefore, $E$ is complete and the calculation also proves closedness of $S_{\lambda,B}$.   
\end{proof}

\begin{remark}
As $S_{\lambda,B} f\in \ker(\At^*-\lambda)$, we have $\norm{S_{\lambda,B} f}_F^2=(1+|\lambda|^2)\norm{S_{\lambda,B} f}_H^2$, so 
\[ \trinorm{f}_E^2:=\norm{f}_\cH^2+\norm{S_{\lambda,B} f}_H^2 \] gives an equivalent norm on $E$.
\end{remark}

\begin{corollary}
If $\Ran (\Gamma_1-B\Gamma_2)=\cH$, then $S_{\lambda,B}:\cH\to  D(\At^*)$ is continuous. In particular, $S_{\lambda, 0}$ is continuous.
\end{corollary}

\begin{proof}
This follows from the Closed Graph Theorem. See for example \cite[Theorem 4.2-I]{Taylor}.
\end{proof}

For the case $\Ran (\Gamma_1-B\Gamma_2)=\cH$, we now want to give a representation of the adjoint of $S_{\lambda,B}$. 
We start with an abstract result:

\begin{lemma}\label{restriction}
Let $M_0\subseteq M$ be a closed subspace of the Hilbert space $M$ and let $N$ be another Hilbert space. 
Suppose $T_1:M_0\to N$ is invertible and $T_2:M\to N$ is such that
\[ (f,h)_M = (f,T_1^{-1}T_2h)_M  \quad \hbox{ for all } f\in M_0,\ h\in M.\] 
Then $T_1=T_2\vert_{M_0}$.
\end{lemma}

\begin{proof}
Let $M=M_0\oplus M_0^\perp$ and $P:M\to M_0$ the orthogonal projection. Then we have $P=T_1^{-1}T_2$ or $T_1P=T_2$ on $M$. 
Therefore, $T_1=T_2$ on $M_0$.  
\end{proof}

\begin{theorem}
Assume $\rho(A_B)\neq\emptyset$. 
Then $A_B^*=\At_{B^*}$.
\end{theorem}

\begin{proof}
Let $u\in D(A_B)$, $v\in D(\At_{B^*})$. Then (\ref{Green}) implies
\begin{equation*}
  (A_B u, v)_H - (u,\At_{B^*}v)_H = (\Gamma_1 u, \Gammat_2 v)_\cH - (\Gamma_2 u, \Gammat_1v)_\cK =
 (B\Gamma_2 u, \Gammat_2 v)_\cH - (\Gamma_2 u, B^*\Gammat_2v)_\cK=0,
\end{equation*}
so $\At_{B^*}\subseteq A_B^*$. On the other hand, let  $v\in D\left(A_B^*\right)$. 
We need to show $(\Gammat_1-B^*\Gammat_2)v=0$. Let $\lambda\in\rho(A_B)$ and $u=(A_B-\lambda)^{-1}w$
for $w\in D(\At^*)$. Then
\begin{eqnarray*}
  0  = (A_Bu,v)- (u, A_B^*v)
    &=& (\Gamma_1 u, \Gammat_2 v)_\cH - (\Gamma_2 u, \Gammat_1v)_\cK \\
    &=& (B\Gamma_2 u, \Gammat_2 v)_\cH - (\Gamma_2 u, \Gammat_1v)_\cK \\
    &=& (\Gamma_2 u,(B^* \Gammat_2-\Gammat_1) v)_\cK\\
    &=& ((A_B-\lambda)^{-1}w,\Gamma_2^*(B^* \Gammat_2-\Gammat_1)v )_{D(\At^*)}\\
    &=& \left(w,\left((A_B-\lambda)^{-1}\right)^*\Gamma_2^*(B^* \Gammat_2-\Gammat_1)v \right)_{D(\At^*)},
\end{eqnarray*}
so $\left((A_B-\lambda)^{-1}\right)^*\Gamma_2^*(B^* \Gammat_2-\Gammat_1)v =0$. Since the adjoint of the 
resolvent is the resolvent of the adjoint,  
$\Gamma_2^*(B^* \Gammat_2-\Gammat_1)v =0$.
Surjectivity of $\Gamma_2$ then gives the result.
\end{proof}

\begin{proposition}
 Assume $\Ran (\Gamma_1-B\Gamma_2)=\cH$  and let $\lambda\in\rho(A_B)$. Then the adjoint of
$S_{\lambda,B}$ is given by $S_{\lambda,B}^*:F\to \cH$,
\begin{eqnarray}
 S_{\lambda,B}^*&=&(1+|\lambda|^2)\ \Gammat_2\ (\At_{B^*}-\overline{\lambda})^{-1}.
\end{eqnarray}
\end{proposition}

\begin{proof}
Choose $v\in\ker(\Gammat_1-B^*\Gammat_2)$, $u\in\ker(\At^*-\lambda)$. Then by (\ref{Green}),
 \begin{eqnarray*}
     -  \left(u,(\At_{B^*}-\overline{\lambda})v\right)_H \!\!\! = 
 \left(\At^* u, v\right)_H \!\!\! - \left(u,\At_{B^*}v\right)_H \!\!\!  = 
\left(\Gamma_1 u, \Gammat_2 v\right)_\cH\!\!\! - \left(\Gamma_2 u, B^*\Gammat_2 v\right)_\cK\!\!\!   \\
=\left((\Gamma_1-B\Gamma_2) u, \Gammat_2 v\right)_\cH.
 \end{eqnarray*} 
 As $S_{\lambda,B}:\cH\to F$ is continuous and continuously invertible, both $S_{\lambda,B}^{*}:F\to \cH$ and $(S_{\lambda,B}^{-1})^*:\cH\to F$
 exist and $(S_{\lambda,B}^*)^{-1}=(S_{\lambda,B}^{-1})^*\in\cL (\cH,F)$.
 Let $w=(\At_{B^*}-\overline{\lambda})v$. Since $\lambda\in\overline{\rho(\At_{B^*})}=\rho(A_B)$, $w\in H$ is arbitrary. 
 Now, by the above calculation, 
 \begin{eqnarray}\label{slamadj}
       - \left(u,w\right)_H
 &=& \left((\Gamma_1-B\Gamma_2)\vert_{\ker(\At^*-\lambda)}u,\Gammat_2(\At_{B^*}-\overline{\lambda})^{-1}w\right)_\cH\\ \nonumber
 & = & \left(S_{\lambda,B}^{-1} u, \Gammat_2(\At_{B^*}-\overline{\lambda})^{-1} w\right)_\cH\\ \nonumber
 & = & \left(u, (S_{\lambda,B}^{-1})^* \Gammat_2(\At_{B^*}-\overline{\lambda})^{-1} w\right)_{F} \\ \nonumber
 & = & \left(u, (S_{\lambda,B}^{*})^{-1} \Gammat_2(\At_{B^*}-\overline{\lambda})^{-1} w\right)_{F} \\ \nonumber
 & = & \left(u, (S_{\lambda,B}^{*})^{-1} \Gammat_2(\At_{B^*}-\overline{\lambda})^{-1} w\right)_{H} 
       +  \left(\At^*u,\At^* (S_{\lambda,B}^{*})^{-1} \Gammat_2(\At_{B^*}-\overline{\lambda})^{-1} w\right)_{H}\\ \nonumber
 & = & (1+|\lambda|^2)\left(u, (S_{\lambda,B}^{*})^{-1}\Gammat_2(\At_{B^*}-\overline{\lambda})^{-1} w\right)_{H}.
 \end{eqnarray}
Therefore, by Lemma \ref{restriction}, we have 
$S_{\lambda,B}^*=(1+|\lambda|^2)\ \Gammat_2(\At_{B^*}-\overline{\lambda})^{-1}$.
\end{proof}

\begin{remark}
\begin{enumerate}
\item The factor $(1+|\lambda|^2)$ is somewhat artificial and comes from the choice of the norm in $F$.
\item Note that since (\ref{slamadj}) only holds for  $u\in\ker(\At^*-\lambda)$, $S_{\lambda,B}^*$ is not defined on the whole of 
$D(\At^*)$. Obviously the operator 
\[ T:=(1+|\lambda|^2)\  \Gammat_2\ (\At_{B^*}-\overline{\lambda})^{-1}P\vert_{\ker(\At^*-\lambda)}\]
is a continuous extension of $S_{\lambda,B}^*$ to $D(\At^*)$ and $T^*= P\vert_{\ker(\At^*-\lambda)}^*S_{\lambda,B}$. Here, $P\vert_{\ker(\At^*-\lambda)}$ denotes the orthogonal projection from $H$ onto $\ker(\At^*-\lambda)$.
\end{enumerate}
\end{remark}

\section{Isolated eigenvalues and poles of the $M$-function}\label{sec:Krein}

For a number of results in what follows we will require an abstract unique continuation hypothesis. We say that the operator $\At^*-\lambda$ satisfies the unique continuation hypothesis if
\[\ker(\At^*-\lambda)\cap\ker(\Gamma_1)\cap\ker(\Gamma_2)=\{0\}.\]
Similarly, $A^*-\lambda$ satisfies the unique continuation hypothesis if
\[\ker(A^*-\lambda)\cap\ker(\Gammat_1)\cap\ker(\Gammat_2)=\{0\}.\]
Whenever either of these conditions is required, it will be stated explicitly.

\begin{remark}
Note that these assumptions are satisfied in the PDE case under fairly general conditions on the operator and the domain (\cf for example \cite[Chapter 4]{Miz73}).
\end{remark}

\begin{lemma}\label{UC}
Assume the unique continuation hypothesis holds for $A^*-\overline{\lambda}$. Then the range of $\At^*-\lambda$ is dense in $H$.
\end{lemma}

\begin{proof}
Suppose there exists $\psi\in H$ such that $\langle \psi,(\At^*-\lambda)u\rangle =0$ 
for all $u\in D(\At^*)$. This implies $\psi\in D(\At^{**})=D(\At)$ and 
$(\At-\overline{\lambda})\psi=0$. The unique continuation hypothesis together with (\ref{domains}) implies $\psi=0$.
\end{proof}

The following definition and Laurent series expansion of the resolvent are standard and can be found in \cite{Kato}. They will be required in a later proof. 

\begin{proposition}\label{notation}
Let $T$ be a closed operator on a Banach space $X$, $\lambda$ an isolated point in the spectrum of $T$ and $\Gamma'$ be a closed path in the resolvent set of $T$ separating $\lambda$ from the rest of the spectrum. 
The spectral projection associated with $\lambda$ is defined by
\begin{eqnarray}\label{P}
	P=\frac{1}{2\pi i}\int_{\Gamma'} R(\zeta,T)d\zeta.
\end{eqnarray}
We also define the eigennilpotent associated with $\lambda$
\begin{eqnarray}\label{D}
	D=(T-\lambda)P=\frac{1}{2\pi i}\int_{\Gamma'} (\zeta-\lambda)R(\zeta,T)d\zeta,
\end{eqnarray}
and
\begin{eqnarray}\label{S}
	S=\frac{1}{2\pi i}\int_{\Gamma'}\frac{1}{\zeta-\lambda}R(\zeta,T)d\zeta.
\end{eqnarray}
For $\zeta$ in a neighbourhood of $\lambda$ the Laurent series expansion of the resolvent is given by
\begin{eqnarray}\label{resrep}
	R(\zeta,T)=\frac{P}{\zeta-\lambda}+\sum_{n=1}^\infty \frac{D^n}{(\zeta-\lambda)^{n+1}}-\sum_{n=0}^\infty (\zeta-\lambda)^nS^{n+1}.
\end{eqnarray}
\end{proposition}

Our aim is now to determine the relationship between the behaviour of the $M$-function $M_B$ as an analytic function and isolated eigenvalues of the operator $A_B$. 

\begin{theorem}\label{Laurentseries}
Let $\mu\in\C$ be an isolated eigenvalue of finite algebraic multiplicity of the operator $A_B$. Assume the unique continuation hypothesis holds for $\At^*-\mu$ and $A^*-\overline{\mu}$.
Then $\mu$ is a pole of finite multiplicity of $M_B(\cdot)$
and the order of the pole of $R(\cdot,A_B)$ at $\mu$ is the same as the order of the pole of $M_B(\cdot)$ at $\mu$.
\end{theorem}

\begin{proof}
We use the following representation of the $M$-function using the resolvent:
\begin{equation}\label{mfn_w}
	M_B(\lambda)f =\Gamma_2\left( I-(A_B-\lambda)^{-1}(\At^*-\lambda)\right)w,
\end{equation}
where $w$ is any element in $D(\At^*)$ such that $(\Gamma_1-B\Gamma_2)w=f$. Obviously,
 any pole of the $M$-function has to be a pole of at least the same order of the resolvent.
It remains to show that the order of the singularity of the pole of the resolvent is preserved despite the presence of the other operators on the right hand side. To do this, we look at the Laurent series expansion. 

Let $\mu$ be an isolated eigenvalue of finite algebraic multiplicity of the operator $A_B$. In this case, there exists $m$ such that the resolvent has a pole of order $m+1$ at $\mu$ and, using the notation from Proposition \ref{notation}, for $\lambda$ in a neighbourhood of $\mu$ the singular part in the representation of the resolvent (\ref{resrep}) is given by  
\[\sum_{n=1}^m \frac{D^n}{(\lambda-\mu)^{n+1}}.\]
In particular, $D^{m+1}=0$ and $D^m\neq 0$. Therefore, there exists $\widetilde{\varphi}$ such that $D^m\widetilde{\varphi}\neq 0$ and $D^m\widetilde{\varphi}$ solves 
\begin{eqnarray*}
  	\left\{ \begin{array}{r@{\; =\;}ll} 
    (\At^*-\mu) u&0 & \\
    (\Gamma_1-B\Gamma_2) u&0& 
\end{array} \right.
\end{eqnarray*}
\ie $D^m\widetilde{\varphi}$ is an eigenfunction of $A_B$ with eigenvalue $\mu$.
We want to show that after substituting the expansion of the resolvent (\ref{resrep}) into $M_B(\mu)$, the most singular term is non-trivial, \ie $\Gamma_2 D^m (\At^*-\mu)w\neq 0$ for some $w\in D(\At^*)$.

First, we show that there exists $\varphi\in H$ satisfying $D^m\varphi\neq 0$ such that the problem $(\At^*-\mu)u=\varphi$ is solvable and $(\Gamma_1-B\Gamma_2)u\neq 0$.
To see this, choose $\widetilde{\varphi}$ such that $D^m\widetilde{\varphi}\neq 0$ and approximate it by $(\varphi_n)\subset \rg (\At^*-\mu)$ which is possible by Lemma \ref{UC}. Since $D^m:H\to H$ is continuous, $D^m\varphi_n\to D^m\widetilde{\varphi}$ and for $N$ sufficiently large, $D^m\varphi_N\neq 0$. Simply choose $\varphi=\varphi_N$.
Now assume $u$ solves $(\At^*-\mu)u=\varphi$ and $(\Gamma_1-B\Gamma_2) u=0$. Then $u\in D(A_B)$ and 
\[ 0=D^{m+1}u=D^m(A_B-\mu)u=D^m\varphi\neq 0,\]
giving a contradiction.

Now we can choose $w$  in (\ref{mfn_w}) as the solution $u$ we have just found. Then $M_B(\lambda)(\Gamma_1-B\Gamma_2) u$ contains the term 
\[\frac{\Gamma_2 D^m (\At^*-\lambda)u}{(\lambda-\mu)^{m+1}}=\frac{\Gamma_2 D^m \left((\At^*-\mu)u-(\lambda-\mu)u\right)}{(\lambda-\mu)^{m+1}},\]
so the most singular term in is of order $(\lambda-\mu)^{-m-1}$ and given by
\[(\lambda-\mu)^{-m-1}\Gamma_2 D^m (\At^*-\mu)u=(\lambda-\mu)^{-m-1}\Gamma_2D^m\varphi.\]
Now $D^m\varphi$ is a (non-trivial) eigenfunction of $A_B$ so by the unique continuation hypothesis, $\Gamma_2D^m\varphi\neq 0$.
\end{proof}

Under slightly stronger hypotheses, we will show next that 
isolated eigenvalues of $A_B$ correspond precisely to isolated poles of the $M$-function. 
We start by proving some identities involving the $M$-function.  
For the $M$-functions associated with two different boundary conditions we have the following identity:
\begin{proposition}
For $\lambda\in\rho(A_B)\cap\rho(A_{B+C})$, we have on $\Ran(\Gamma_1-B\Gamma_2)$
\begin{equation}\label{mfn2}
M_{B+C}(\lambda)(I-CM_B(\lambda))=M_B(\lambda).
\end{equation}
Correspondingly, we have
\begin{equation}\label{slam2}
S_{\lambda,B+C}(I-C\Gamma_2 S_{\lambda,B})=S_{\lambda,B}\quad \hbox{ on }\quad \Ran(\Gamma_1-B\Gamma_2).
\end{equation}
\end{proposition}

\begin{proof}
We prove (\ref{slam2}). Then (\ref{mfn2}) follows by applying $\Gamma_2$ to both sides. Let $f\in\Ran(\Gamma_1-B\Gamma_2)$, then
$(\Gamma_1-B\Gamma_2)S_{\lambda,B}f=f$, so
\[S_{\lambda,B+C}(I-C\Gamma_2 S_{\lambda,B})f=S_{\lambda,B+C}(\Gamma_1-B\Gamma_2-C\Gamma_2)S_{\lambda,B}f= S_{\lambda,B}f,\]
since $S_{\lambda,B}f\in\ker(\At^*-\lambda)$.
\end{proof}

The next proposition gives a representation of the $M$-function in terms of the resolvent.
\begin{proposition}
Let $\lambda,\lambda_0\in\rho(A_B)$. Then on $\Ran(\Gamma_1-B\Gamma_2)$
\begin{eqnarray}\label{mfn}
M_B(\lambda) &=& \Gamma_2\left( I+(\lambda-\lambda_0)(A_B-\lambda)^{-1}\right)S_{\lambda_0,B}\ 
=\ \Gamma_2(A_B-\lambda_0)(A_B-\lambda)^{-1}S_{\lambda_0,B}.
\end{eqnarray}
\end{proposition}

\begin{proof}
Given $f\in\Ran(\Gamma_1-B\Gamma_2)$, choose $u\in D(\At^*)$ such that $(\Gamma_1-B\Gamma_2)u=f$.
Set
 \[w=u-(A_B-\lambda)^{-1}(\At^*-\lambda)u.\] 
Then $w\in\ker(\At^*-\lambda)$, $(\Gamma_1-B\Gamma_2)w=f$ and $w$
is the unique function with these properties, as $\lambda\in\rho(A_B)$. Moreover, $M_B(\lambda)f=\Gamma_2 w$.  Choose $u=S_{\lambda_0,B}f$.
Then
\begin{eqnarray*}   
M_B(\lambda)f &=& \Gamma_2\left( I-(A_B-\lambda)^{-1}(\At^*-\lambda)\right)S_{\lambda_0,B}f\\ 
             &=& \Gamma_2\left( I+(\lambda-\lambda_0)(A_B-\lambda)^{-1}\right)S_{\lambda_0,B}f\\ 
          &=& \Gamma_2(A_B-\lambda_0)(A_B-\lambda)^{-1}S_{\lambda_0,B}f. 
\end{eqnarray*}
\end{proof}

We now give a representation of the resolvent in terms of the $M$-function. This type of formulae are usually called Kre\u{\i}n's formulae.

\begin{theorem}
Let $B,C\in \cL(\cK,\cH)$, $\lambda\in\rho(A_B)\cap\rho(A_C)\cap\rho(A_{B+C})$. Then 
\begin{eqnarray}\label{Krein}
(A_B-\lambda)^{-1} &=& (A_C-\lambda)^{-1}-S_{\lambda,B+C}(I-CM_B(\lambda))(\Gamma_1-B\Gamma_2)(A_C-\lambda)^{-1}\\ \nonumber
                   &=& (A_C-\lambda)^{-1}-S_{\lambda,B+C}(I-CM_B(\lambda))(C-B)\Gamma_2(A_C-\lambda)^{-1}.
\end{eqnarray}
\end{theorem}

\begin{proof}
Let $u\in H$. Set $v:=\left((A_B-\lambda)^{-1}-(A_C-\lambda)^{-1}\right)u$. Since $v\in\ker(\At^*-\lambda)$, we have 
$M_B(\lambda)(\Gamma_1-B\Gamma_2)v=\Gamma_2 v$. Then
\begin{eqnarray}\label{Kreinpf}
\left(\Gamma_1-(B+C)\Gamma_2\right)v &=& \left[\Gamma_1-B\Gamma_2-CM_B(\lambda)(\Gamma_1-B\Gamma_2)\right] v\\ \nonumber
                          &=& (I-CM_B(\lambda))(\Gamma_1-B\Gamma_2) v\\ \nonumber
                          &=& (I-CM_B(\lambda))(\Gamma_1-B\Gamma_2) \left((A_B-\lambda)^{-1}-(A_C-\lambda)^{-1}\right)u\\ \nonumber
                          &=& -(I-CM_B(\lambda))(\Gamma_1-B\Gamma_2) (A_C-\lambda)^{-1}u.
\end{eqnarray}
Set $f:=-(I-CM_B(\lambda))(\Gamma_1-B\Gamma_2) (A_C-\lambda)^{-1}u$. Then, by (\ref{Kreinpf}), $f\in\Ran(\Gamma_1-(B+C)\Gamma_2)$
and $S_{\lambda, B+C}f=v=\left((A_B-\lambda)^{-1}-(A_C-\lambda)^{-1}\right)u$. Therefore, 
\begin{eqnarray*}
(A_B-\lambda)^{-1}= (A_C-\lambda)^{-1}-S_{\lambda, B+C}(I-CM_B(\lambda))(\Gamma_1-B\Gamma_2) (A_C-\lambda)^{-1}.	
\end{eqnarray*}
\end{proof}

\begin{remark}
If $\lambda\in\rho(A_B)\cap\rho(A_C)\cap\rho(A_{B-C})$, then we have 
\[(A_B-\lambda)^{-1} = (A_C-\lambda)^{-1}-S_{\lambda,B-C}(I+CM_B(\lambda))(C-B)\Gamma_2(A_C-\lambda)^{-1}. \]
\end{remark}

The case $B=0$ is particularly simple:
\begin{corollary}
Let $C\in \cL(\cK,\cH)$, $\lambda\in\rho(A_0)\cap\rho(A_C)$. Then 
\begin{eqnarray*}
(A_0-\lambda)^{-1} &=& (A_C-\lambda)^{-1}-S_{\lambda,C}(I-CM_0(\lambda))\Gamma_1(A_C-\lambda)^{-1}.
\end{eqnarray*}
\end{corollary}

We our now ready to prove our main result.
\begin{theorem}
Let $\mu\in\C$. We assume that $\rho(A_B)\neq\emptyset$ and that
there exist operators $B,C\in \cL(\cK,\cH)$ such that $\mu\in\rho(A_C)\cap\rho(A_{B+C})$ or $\mu\in\rho(A_C)\cap\rho(A_{B-C})$.
Then $\mu$ is an isolated eigenvalue of finite algebraic multiplicity of the operator $A_B$ if and only if
$\mu$ is a pole of finite multiplicity of $M_B(\cdot)$.
In this case, the order of the pole of $R(\cdot,A_B)$ at $\mu$ is the same as the order of the pole of $M_B(\cdot)$ at $\mu$.
\end{theorem}

\begin{proof}
Let $\mu$ be an isolated eigenvalue of finite algebraic multiplicity $m$ of the operator $A_B$. 
Then, since  $\mu\in\rho(A_C)\cap\rho(A_{B\pm C})$, and $S_{\lambda,B\pm C}$ is analytic in $\lambda$ by Proposition \ref{slamanaprop},
 (\ref{Krein}) implies that $M_B(\cdot)$ must 
have a pole of at least order $m$ at $\mu$, while (\ref{mfn}) implies that the pole is at most of order $m$.

Similarly, if $M_B(\cdot)$ has a pole of order $m$ at $\mu$, (\ref{mfn}) implies that the resolvent of $A_B$ must have a pole of order 
at least $m$ at $\mu$, while (\ref{Krein}) implies that the pole is at most of order $m$. Therefore, $\mu$ is an eigenvalue of $A_B$ 
(\cf for example \cite[Section 3.6.5]{Kato}). 
\end{proof}

\begin{remark}
Note that the assumption that $C$ can be chosen such that  $\mu\in\rho(A_C)$ implies the unique continuation property for $\At^*-\mu$.

To see this, let $u\in\ker(\At^*-\mu)\cap\ker(\Gamma_1)\cap\ker(\Gamma_2)$. 
Then $u\in\ker(\Gamma_1-C\Gamma_2)$, so $u\in D(A_C)$ and $(A_C-\mu)u=0$, so $u=(A_C-\mu)^{-1}(A_C-\mu)u=0$.
\end{remark}

\section{Behaviour of the $M$-function near the essential spectrum}

By the essential spectrum of an operator $\sigma_{ess}$, we denote all points in the spectrum that are not isolated eigenvalues of finite multiplicity. In this section we will investigate what can be said about the essential spectrum from the behaviour of the $M$-function. In the case of symmetric operators, these questions have been addressed by Brasche, Malamud and Neidhardt in \cite{Brasche}.


\begin{theorem}\label{cutlinethm}
Let $k\in\C$ such that there exists $\eps_0>0$, with $k\pm i\eps\in\rho(A_B)$ for all $0<\eps<\eps_0$.
Suppose there is a linear subspace $\fH\subseteq H$ such that $\fH\cap D(A^*)$ is dense
in $H$ and  
\begin{enumerate}
	\item for every $f\in\Ran(\Gamma_1-B\Gamma_2)$ we can find $F\in \fH\cap D(\At^*)$ satisfying
\begin{itemize}
	\item $(\Gamma_1-B\Gamma_2)F=f$,
	\item $u:=(\At^*-k)F\in \fH$;
\end{itemize}
        \item $\left. (\Gammat_1-B^*\Gammat_2)\right|_{\fH\cap D(A^*)}$ is surjective;
        \item for all $v\in\fH\cap D(A^*)$, $A^*v\in \fH$;
	\item $\lim_{\eps\to 0}((A_B-(k\pm i\eps))^{-1}w,v) \hbox{ exists  for all } w, v \in \fH.$
\end{enumerate}
Then for all $f\in\Ran(\Gamma_1-B\Gamma_2)$ the weak limits $M_B(k\pm i0)f:=\wlim_{\eps\to 0}M_B(k\pm i\eps)f$ exist. 
Moreover, \[(A_B-(k+ i0))^{-1}u=(A_B-(k- i0))^{-1}u \quad \hbox{ implies }\quad M_B(k+ i0)f=M_B(k- i0)f.\]
Here, the left hand equality  is to be interpreted as 
\[ \lim_{\eps\to 0} \left((A_B-(k+ i\eps))^{-1}u,v\right) = \lim_{\eps\to 0} \left((A_B-(k- i\eps))^{-1}u,v\right) \quad\hbox{ for all } v\in\fH.\]
\end{theorem}

\begin{remark}
In the case of an elliptic PDE in an unbounded domain with finite boundary, the subspace $\fH$ could consist of all finitely 
supported functions.  

Condition (4) is our main assumption, while (1) is a kind of inverse trace theorem and (2) and (3) are technical assumptions. 
\end{remark}

\begin{proof}
Given $f\in\Ran(\Gamma_1-B\Gamma_2)$, choose $F\in \fH$ such that $(\Gamma_1-B\Gamma_2)F=f$. Set \[w_{\eps,\pm}:=F-(A_B-(k\pm i\eps))^{-1}(\At^*-(k\pm i\eps))F.\] 
Then $w_{\eps,\pm}\in\ker(\At^*-(k\pm i\eps))$, $M_B(k\pm i\eps)f=\Gamma_2 w_{\eps,\pm}$ and $\Gamma_1 w_{\eps,\pm}=(\Gamma_1-B\Gamma_2+B\Gamma_2) w_{\eps,\pm}=(I+BM_B(k\pm i\eps))f$.
Green's identity (\ref{Green}) for any $v\in D(A^*)$ gives
\begin{eqnarray*}
-\left(w_{\eps,\pm},(A^*-(\overline{k}\mp i\eps))v\right)_H
&=&\left((\At^*-(k\pm i\eps))w_{\eps,\pm}, v\right)_H -\left( w_{\eps,\pm},(A^*-(\overline{k}\mp i\eps))v\right)_H \\ \nonumber
&=&\left(\Gamma_1 w_{\eps,\pm},\Gammat_2 v\right)_\cH - \left(\Gamma_2 w_{\eps,\pm},\Gammat_1 v\right)_\cK  \\ \nonumber
&=&\left((I+BM_B(k\pm i\eps))f,\Gammat_2 v\right)_\cH - \left(M_B(k\pm i\eps)f,\Gammat_1 v\right)_\cK  \\ \nonumber
&=&\left(f,\Gammat_2 v\right)_\cH - \left(M_B(k\pm i\eps)f,(\Gammat_1 -B^*\Gammat_2)v\right)_\cK .
\end{eqnarray*}
Setting $u=(\At^*-k)F$ and inserting our expression for $w_{\eps,\pm}$ on the left hand side, the equation becomes
\begin{eqnarray}\label{Greeneps}
\quad \left(F-(A_B-(k\pm i\eps))^{-1}(u\mp i\eps F),(A^*-(\overline{k}\mp i\eps))v\right)_H
&=& - \left(f,\Gammat_2 v\right)_\cH \\ \nonumber && + \left(M_B(k\pm i\eps)f,(\Gammat_1 -B^*\Gammat_2)v\right)_\cK .
\end{eqnarray}

Now assume $v\in \fH\cap D(A^*)$. Since $u, F\in \fH$, we can take limits on the left hand side. The assumption that
$\left. (\Gammat_1 -B^*\Gammat_2)\right|_{\fH\cap D(A^*)}$ is surjective then gives weak convergence 
of $M_B(k\pm i\eps)f $ in $\cK$ and we get
\begin{eqnarray}\label{Greenlim}
\left(F-(A_B-(k\pm i0))^{-1}u,(A^*-\overline{k})v\right)_H 
&=& - \left(f,\Gammat_2 v\right)_\cH  + \left(M_B(k\pm i0)f,(\Gammat_1 -B^*\Gammat_2)v\right)_\cK .
\end{eqnarray}
Furthermore,
\begin{eqnarray}\label{Greenlimdiff}
\left(((A_B-(k+ i0))^{-1} - (A_B-(k- i0))^{-1})u,(A^*-\overline{k})v\right)_H \\ \nonumber
&\hspace{-60mm}=&\hspace{-28mm} - \left((M_B(k+ i0) -M_B(k- i0))f,(\Gammat_1 -B^*\Gammat_2)v\right)_\cK .
\end{eqnarray}
Since
$\left. (\Gammat_1 -B^*\Gammat_2)\right|_{\fH\cap D(A^*)}$ is surjective, equality of the weak limits of the resolvent implies equality of the weak limits of the $M$-function. 
\end{proof}

We would like to prove a converse of Theorem \ref{cutlinethm}, \ie determine the behaviour of the resolvent from that of the $M$-function. However, we only get the following partial results:
\begin{proposition} Assume the unique continuation hypothesis holds for $\At^*-k$ and $A^*-\overline{k}$ and that the weak limits \[M_B(k\pm i0)g:=\wlim_{\eps\to 0}M_B(k\pm i\eps)g\]
exist for every $g\in\Ran(\Gamma_1-B\Gamma_2)$ and that there exists some $f\in\Ran(\Gamma_1-B\Gamma_2)$ such that 
\[M_B(k+ i0)f\neq M_B(k - i0)f.\] 
Then  $k\in\sigma_{ess}(A_B)$. 
\end{proposition}

\begin{remark}
Note that in \cite{Brasche} it is shown that for symmetric operators $\Im (M_B(k+ i0)f,f)>0$ implies $k\in\sigma_{ess}(A_B)$.
\end{remark}

\begin{proof}
As in the proof of Theorem \ref{cutlinethm}, we arrive at equation (\ref{Greeneps}). By assumption, the limit on the right hand side exists.
Assume that $k\in\rho(A_B)$. Then we can take limits on the left hand side and get equation (\ref{Greenlimdiff}) with the l.h.s. equal to 0 
contradicting   $M_B(k+ i0)f\neq M_B(k - i0)f$. Thus $k\in\sigma(A_B)$ and $k$ is not in the isolated point spectrum, as the weak limits of the $M$-function exist which would contradict Theorem \ref{Laurentseries}.
\end{proof}

\begin{remark}
 The problem in getting a stronger statement lies in the fact that the
 $M$-function does not contain all the singularities of the resolvent, but
 only those that are contained on a certain subspace. We plan to  
 discuss this
 topic and other properties related to the continuous spectrum and
 behaviour of the $M$-function in a forthcoming paper.
%
%
%
\end{remark}

In what follows, we will show that for a block operator matrix it is possible
to have a dense proper subspace $\fH$ for which the weak limit of the $M$-functions exists, 
but the weak limit for the resolvents does not exist. We also hope that this example, 
demonstrating the calculation of the $M$-function in a non-trivial block 
operators matrix setting,  is of independent interest.

\subsection*{A block matrix-differential operator related to the Hain-L\"{u}st operator}

Let 
\begin{equation}
\At^* = \left(\begin{array}{cc} -\frac{d^2}{dx^2}+q(x) & w(x) \vspace{2pt}\\
 w(x) & u(x) \end{array}\right), 
\label{eq:hl1} 
\end{equation}
where $q$, $u$ and $w$ are $L^\infty$-functions, and the domain of
the operator is given by
\begin{equation}
D(\At^*) = H^2(0,1)\times L^2(0,1). 
\label{eq:hl2} 
\end{equation}
Also let 
\begin{equation}
 A^* = \left(\begin{array}{cc} -\frac{d^2}{dx^2}+\overline{q(x)} & \overline{w(x)} \vspace{2pt}\\ 
 \overline{w(x)} & \overline{u(x)} \end{array}\right).
\label{eq:hl3}
\end{equation}
It is then easy to see that
\begin{eqnarray}
\llangle \At^*\left(\begin{array}{c} y \\ z \end{array}\right),
 \left(\begin{array}{c} f \\ g \end{array}\right)\rrangle
 -  \llangle \left(\begin{array}{c} y \\ z \end{array}\right),
  A^*\left(\begin{array}{c} f \\ g \end{array}\right)\rrangle\nonumber & & \\
 & \hspace{-9cm} = & \hspace{-4.5cm}\llangle \Gamma_1\left(\begin{array}{c} y \\ z \end{array}\right),
 \Gamma_2\left(\begin{array}{c} f \\ g \end{array}\right)\rrangle
 - \llangle \Gamma_2\left(\begin{array}{c} y \\ z \end{array}\right),
 \Gamma_1\left(\begin{array}{c} f \\ g \end{array}\right)\rrangle, 
\label{eq:hl4}
\end{eqnarray}
where 
\[ \Gamma_1\left(\begin{array}{c} y \\ z \end{array}\right)
  = \left(\begin{array}{c} -y'(1) \\ y'(0) \end{array}\right),
\;\;\;
 \Gamma_2\left(\begin{array}{c} y \\ z \end{array}\right)
  = \left(\begin{array}{c} y(1) \\ y(0) \end{array}\right).
\]
Consider the operator
\begin{equation}
 A_{\alpha\beta} := \left. \At^*\right|_{\mbox{ker}(\Gamma_1-B\Gamma_2)},
\label{eq:hl5}
\end{equation}
where $B = \left(\begin{array}{cc} \cot\beta & 0 \\ 
                0 & -\cot\alpha \end{array}\right)$.
It is known (see, e.g., \cite{Langer}) that 
$\sigma_{ess}(A_{\alpha\beta}) = \mbox{essran}(u)$. This result is independent
of the choice of boundary conditions.

We now calculate the function $M(\lambda)$ such that
\[ M(\lambda)(\Gamma_1 - B \Gamma_2)\left(\begin{array}{c} y \\ z 
 \end{array}\right) = \Gamma_2\left(\begin{array}{c} y \\ z 
 \end{array}\right) \]
for $ \left(\begin{array}{c} y \\ z  \end{array}\right)\in
 \mbox{ker}(\At^*-\lambda)$. In our calculation we assume
that $\lambda\not\in \sigma_{ess}(A_{\alpha\beta})$. The condition
$ \left(\begin{array}{c} y \\ z  \end{array}\right)\in
 \mbox{ker}(\At^*-\lambda)$ yields the equations
\[ -y'' + (q-\lambda)y + wz = 0;\;\;\
 wy + (u-\lambda)z = 0 \]
which, in particular, give
\begin{equation}
 -y'' + (q-\lambda)y + \frac{w^2}{\lambda-u}y = 0.
\label{eq:hl10}
\end{equation}
The linear space $\mbox{ker}(\At^*-\lambda)$ is therefore
spanned by the functions 
$ \left(\begin{array}{c} y_1 \\ wy_1/(\lambda-u)  \end{array}\right)$
and 
$ \left(\begin{array}{c} y_2 \\ wy_2/(\lambda-u)  \end{array}\right)$
where $y_1$ and $y_2$ are solutions of the initial value problems
consisting of the differential equation (\ref{eq:hl10}) equipped
with initial conditions
\begin{equation}
 y_1(0) = \cos\alpha, \;\; y_1'(0) = \sin\alpha, \label{eq:hl10a}
\end{equation}
\begin{equation}
 y_2(0) = -\sin\alpha, \;\; y_2'(0) = \cos\alpha. \label{eq:hl10b}
\end{equation}
A straightforward calculation shows that
\[ \left(\begin{array}{c} y(1) \\ y(0) \end{array}\right)
 = \left(\begin{array}{cc} m_{11}(\lambda) & m_{12}(\lambda) \\
 m_{21}(\lambda) & m_{22}(\lambda) \end{array}\right)
\left(\begin{array}{c} -y'(1) -\cos\beta\ y(1)/\sin\beta\\ y'(0)
 + \cos\alpha\ y(0)/\sin\alpha \end{array}\right). \]
Note that the $y_j$ depend on $x$ and $\lambda$ but that the
$\lambda$-dependence is suppressed in the notation, except when necessary. 
Another elementary calculation now shows that
\begin{equation}
 m_{11}(\lambda) = -\frac{y_2(1,\lambda)}{y_2'(1,\lambda) + 
 \cot\beta\ y_2(1,\lambda)}, 
\label{eq:hl7} 
\end{equation}
\begin{equation}
 m_{21}(\lambda) = m_{12}(\lambda) = \frac{\sin\alpha}{y_2'(1,\lambda) + 
 \cot\beta\ y_2(1,\lambda)}, 
\label{eq:hl8} 
\end{equation}
\begin{equation}
 m_{22}(\lambda) = \sin\alpha\cos\alpha+\sin^2\alpha\left\{
\frac{y_1'(1,\lambda) + \cot\beta\ y_1(1,\lambda)}{y_2'(1,\lambda) + 
 \cot\beta\ y_2(1,\lambda)}\right\}.
\label{eq:hl9} 
\end{equation}

As an aside, notice that all these expressions contain a denominator 
$y_2'(1,\lambda) +  \cot\beta\ y_2(1,\lambda)$ and that 
$\lambda\not\in \mbox{essran}(u)$ is an eigenvalue precisely when
this denominator is zero. 

We now fix $k\in \mbox{essran}(u)$, let $\lambda = k \pm i\eps$, 
and consider the limits $\lim_{\eps\searrow 0}M(k\pm i\eps)$.
For simplicity we consider the case in which $u$ is injective and
$k=u(x_0)$ for some $x_0\in (0,1)$ and we suppose that $w(x)=0$ for
$x \in (x_0-\delta,x_0+\delta)$ for some small $\delta>0$. In 
this case the coefficient 
\[ \frac{w(x)}{u(x)-\lambda} \]
is well defined as a function of $x$ for all $\lambda$ in a 
punctured neighbourhood in $\C$ of the point $k=u(x_0)$: in particular,
$w(x)/(u(x)-\lambda)$ is identically zero for all $\lambda\neq k$,
for all $x\in (x_0-\delta,x_0+\delta)$. Consequently the solutions
$y_1(x,\lambda)$ and $y_2(x,\lambda)$ are well defined for all 
$x\in [0,1]$, for all $\lambda$ in a neighbourhood of $k=u(x_0)$.
The $M$-function may have an isolated pole at some point 
$\lambda$ near $k$ if $y_2'(1,\lambda) +  \cot\beta\ y_2(1,\lambda)$
happens to be zero; such a pole will be an eigenvalue of the operator
$A_{\alpha\beta}$ embedded in the essential spectrum and therefore
a more complicated singularity of $(A_{\alpha\beta}-\lambda)^{-1}$.
Embedded eigenvalues may occur even without the hypothesis that $w$ 
vanishes on some subinterval $(x_0-\delta,x_0+\delta)$: see \cite{BrownLangerMarletta}.
However embedded eigenvalues are atypical and are generally destroyed by
an arbitrarily small perturbation to the problem. In the absence of any 
embedded eigenvalues, $M(\lambda)$ will be analytic in the neighbourhood
$u(x_0-\delta,x_0+\delta)$ of the point $k=u(x_0)$ and we shall have, 
in the sense of norm limits,
\[ \lim_{\eps\searrow 0}M(\mu+i\eps)
 = \lim_{\eps\searrow 0}M(\mu-i\eps) \;\;\;
\forall \mu \in u(x_0-\delta,x_0+\delta). \]
For the resolvent, suppose that 
\[ \left(\begin{array}{c} y \\ z \end{array}\right) 
 = (A_{\alpha\beta}-\lambda)^{-1}\left(\begin{array}{c} f_1 \\ f_2 \end{array}
\right). \]
Then $y$ must satisfy
\[ -y'' + (q-\lambda)y -\frac{w^2}{u-\lambda}y 
 = f_1 - \frac{w}{u-\lambda}f_2, \]
together with the boundary conditions, which is a uniquely solvable problem
in the absence of embedded eigenvalues (recall that $w/(u-\lambda)$ is 
well defined as a function of $x$ for all $\lambda$ in a neighbourhood of
$k$). In particular, $y(x,\lambda)$ does not have a singularity of
any type at $\lambda = u(x_0)$.

Now $z$ is given by
\be z = \frac{f_2}{u-\lambda} - \frac{w}{u-\lambda}y. \label{eq:zed}\ee
We examine the question of existence of weak limits of the type described
in Theorem \ref{cutlinethm}:
\[ \lim_{\eps\searrow 0}\langle (A_B-\lambda)^{-1}f,g\rangle \]
where $f=(f_1,f_2)$ and $g=(g_1,g_2)$ lie in some space $\fH$
and $\lambda = u(x_0)\pm i\eps$. Evidently the first component $y$ 
of the vector $(A_B-\lambda)^{-1}f$ will cause
no problems whatever ${\fH}$ we choose:
\[ \int_0^1 y(x,\lambda)\overline{g_1(x)}dx \]
will be analytic in a neighbourhood of $\lambda=u(x_0)$. Thus we turn to 
the second component $z(x,\lambda)$. Take ${\fH}$ to be the space of 
two-component smooth functions. Suppose that $u$ is differentiable at 
$x_0\in (0,1)$ with $u'(x_0)\neq 0$. If $z$ is given by (\ref{eq:zed}) 
then the inner product
\[ \int_0^1 z(x,\lambda)\overline{g_2(x)}dx \]
with $\lambda = u(x_0)+ i\eps$ has a limit as $\eps$ tends to
zero from above; similarly as it has a (generally different) limit as
$\eps$ tends to zero from below. The difference of the limits is
\be 2\pi i f_2(x_0)g_2(x_0). \label{eq:limdif}\ee
However the $M$-function has no singularity at all. We have therefore 
constructed an example in which the resolvent has non-equal weak limits
but the $M$-function has equal norm limits. 

It is worth emphasizing that for this example,
\[ \overline{\mbox{Ran}(A^*-\overline{k})} = H. \]
This is not enough to avoid the phenomenon that some singularities of the
resolvent are `canceled' in the $M$-function.

\section{Relatively bounded perturbations}

Let $U$ be a symmetric operator in $H$ and $(\cH,\Gamma_1,\Gamma_2)$ be a boundary value space for $U$ (\cf \cite[pp 155]{Gorbachuk}).
Assume that $V$ is symmetric with the following properties:
\begin{itemize}
\item $V$ is relatively $U$-bounded with relative bound less than 1
\item $V^*$ is relatively $U^*$-bounded with relative bound less than 1
\end{itemize}
We will show that in this case it is sufficient to consider boundary operators only associated with the symmetric part $U$ of the operator $A=U+iV$.

\begin{example}
Let $U$ be a symmetric second order elliptic differential operator on a smooth domain $\Omega\subseteq\R^n$ with $D(U)=H^2_0(\Omega)$. 
If $n>1$, only operators of the form $Vu=qu$, $q\in L^\infty(\Omega,\R)$ satisfy these conditions. If $n=1$, $V$ can also involve 
first order terms.
\end{example}

Let $A=U+iV$ and $\At=U-iV$. By the assumptions on $V$, $D(A)=D(\At)=D(U)$ 
and $A^*=U^*-iV^*$, $\At^*=U^*+iV^*$. with $D(A^*)=D(\At^*)=D(U^*)$.
Moreover, $A\subseteq\At^*$ and $\At\subseteq A^*$. For $B\in\cL(\cH)$, let $A_B=\At^*|_{\ker(\Gamma_1-B\Gamma_2)}$ and define 
$M_B(\lambda)$ and $S_{\lambda,B}$ as before with the boundary operators $\Gamma_1, \Gamma_2$ now only associated with the symmetric part of $A$.
Then all the results of Section \ref{sec:Krein} hold in this situation as well and the 
proofs are identical as the specific form of the Green formula plays no role in their derivation. Therefore, we have

\begin{theorem}
Let $\mu\in\C$ be an isolated eigenvalue of finite algebraic multiplicity of the operator $A_B$. Assume the unique continuation hypothesis holds for $\At^*-\mu$ and $A^*-\overline{\mu}$.
Then $\mu$ is a pole of finite multiplicity of $M_B(\cdot)$
and the order of the pole of $R(\cdot,A_B)$ at $\mu$ is the same as the order of the pole of $M_B(\cdot)$ at $\mu$.
\end{theorem}

\begin{proposition}
For $\lambda\in\rho(A_B)\cap\rho(A_{B+C})$, we have
\begin{equation*}
M_{B+C}(\lambda)(I-CM_B(\lambda))=M_B(\lambda).
\end{equation*}
Correspondingly, we have
\begin{equation*}
S_{\lambda,B+C}(I-C\Gamma_2 S_{\lambda,B})=S_{\lambda,B}.
\end{equation*}
\end{proposition}

\begin{proposition}
Let $\lambda,\lambda_0\in\rho(A_B)$. Then 
\begin{eqnarray*}
M_B(\lambda) &=& \Gamma_2\left( I+(\lambda-\lambda_0)(A_B-\lambda)^{-1}\right)S_{\lambda_0,B}\ 
=\ \Gamma_2(A_B-\lambda_0)(A_B-\lambda)^{-1}S_{\lambda_0,B}.
\end{eqnarray*}
\end{proposition}

\begin{proposition}
Let $B,C\in \cL(\cH)$, $\lambda\in\rho(A_B)\cap\rho(A_C)\cap\rho(A_{B+C})$. Then 
\begin{eqnarray*}
(A_B-\lambda)^{-1} &=& (A_C-\lambda)^{-1}-S_{\lambda,B+C}(I-CM_B(\lambda))(\Gamma_1-B\Gamma_2)(A_C-\lambda)^{-1}\\ \nonumber
                   &=& (A_C-\lambda)^{-1}-S_{\lambda,B+C}(I-CM_B(\lambda))(C-B)\Gamma_2(A_C-\lambda)^{-1}.
\end{eqnarray*}
\end{proposition}
\begin{theorem}
Let $\mu\in\C$ and assume there exist operators $B,C\in \cL(\cH)$ such that $\mu\in\rho(A_C)\cap\rho(A_{B+C})$ or 
$\mu\in\rho(A_C)\cap\rho(A_{B-C})$. 
Then $\mu$ is an isolated eigenvalue of finite algebraic multiplicity of the operator $A_B$ if and only if
$\mu$ is a pole of finite multiplicity of $M_B(\cdot)$.
In this case, the order of the pole of $R(\cdot,A_B)$ at $\mu$ is the same as the order of the pole of $M_B(\cdot)$ at $\mu$.
\end{theorem}

\section{Application to PDEs}
The theory previously developed is not immediately applicable to the usual boundary value problems arising in PDEs. The reason is the following:
Consider the case of the Laplacian $A=\Delta$ with $D(A)=H^2_0(\Omega)$ where $\Omega$ is a smooth bounded domain. The usual Green's identity is given by 
\[ \into \left(-\Delta u \overline{v}+u\Delta \overline{v}\right)= \intdo \left(-\frac{\partial u}{\partial \nu} \overline{v} + u  \frac{\partial \overline{v}}{\partial \nu}\right), \quad u,v\in H^2(\Omega). 
\]
However, we want identity (\ref{Green}) to hold for all $u,v\in D(\At^*)=D(A^*)=\{u\in L^2(\Omega): \Delta u \in  L^2(\Omega)\}$ which in general is not even a subset of $H^1(\Omega)$. Therefore, the integral $\intdo \frac{\partial u}{\partial \nu} \overline{v}$ is not well-defined for all these functions. 

The aim of this section is to show that by suitably modifying the boundary operators, our previous results hold for elliptic differential operators of order $2m$. This idea was first used by Vishik \cite{Vis52}. So as not to obscure the ideas with technicalities and notation we will only consider a first order perturbation of the Laplacian. The same method is applicable to any elliptic operator satisfying the conditions given in \cite[\S I.3]{Grubb} by Grubb. In fact, all the results required in the following are taken from that paper. 

Let 
\[
A=\Delta+p\cdot\nabla, \quad D(A)=H^2_0(\Omega), \quad p\in(C^\infty(\overline{\Omega}))^n\]
\[\At=\Delta-\dive (p\ \cdot ), \quad D(\At)=H^2_0(\Omega),\]
where $\Omega$ is a smooth bounded domain. 
Let 
\[
\gamma_1 u= \left[\frac{\partial u}{\partial \nu}+(p\cdot\nu)u\right]\Bigg|_{\partial\Omega}, \quad \gamma_2 u=u\Big|_{\partial\Omega}\]
\[ \gammat_1 v= \frac{\partial v}{\partial \nu}\Big|_{\partial\Omega},\quad \gammat_2 v=v\Big|_{\partial\Omega}\]
Then for $u,v\in H^2(\Omega)$ we have
\[(\At^* u,v)_{L^2(\Omega)} - (u,A^*v)_{L^2(\Omega)} = (\gamma_1 u, \gammat_2 v)_{L^2(\partial\Omega)}- (\gamma_2 u, \gammat_1 v)_{L^2(\partial\Omega)}.\]

It is easy to check that \[ D(\At^*)=\{u\in L^2(\Omega): (\Delta + p\cdot\nabla)u\in L^2(\Omega)\},\] 
\[ D(A^*)=\{v\in L^2(\Omega): \Delta v -\dive  (p\ v)\in L^2(\Omega)\}.\] 

Let $A_D:=\At^*\big|_{\ker \gamma_2}$ be the restriction of $\At^*$ satisfying Dirichlet boundary conditions. 
Similarly, let $\At_D:=A^*\big|_{\ker \gammat_2}$. Then by elliptic regularity, $D(A_D)=H^2(\Omega)\cap H^1_0(\Omega)=D(\At_D)$. 
Without loss of generality, assume that $0\in\rho(A_D)\cap\rho(\At_D)$ (if not, this can be achieved by a translation). By \cite[Lemma II.1.1]{Grubb},
$D(\At^*)=D(A_D)+\ker\At^* $ and $D(A^*)=D(A_D)+\ker A^* $.

\begin{definition}
For $\varphi\in H^{-1/2}(\partial\Omega)$ define $m_0\varphi\in H^{-3/2}(\partial\Omega)$ by 
\[ m_0\varphi=\gamma_1 u=\left(\frac{\partial u}{\partial \nu}+(p\cdot\nu)u\right)\Bigg|_{\partial\Omega},\quad  \hbox{ where $u$ solves } \At^* u= 0,\quad \gamma_2 u=\varphi\]
and let $\mt_0\varphi\in H^{-3/2}(\partial\Omega)$ satisfy
\[ \mt_0\varphi=\gammat_1 v=\frac{\partial v}{\partial \nu}\Big|_{\partial\Omega},\quad \hbox{ where $v$ solves } A^* v= 0,\quad \gammat_2 v=\varphi.\]
\end{definition}

\begin{definition}
For $u\in D(\At^*)$, let \[\Gamma u:= \gamma_1 u-m_0\gamma_2 u\]
and for $v\in D(A^*)$, let \[\Gammat v:= \gammat_1 v-\mt_0\gammat_2 v.\]
\end{definition}

\begin{remark}
\begin{enumerate}
	\item  The operators $m_0, \mt_0, \Gamma$ and $\Gammat$ are well-defined (\cf \cite[\S III.1]{Grubb}).
	\item  $m_0$ and $\mt_0$ are the Dirichlet to Neumann maps associated with $\At^*$ and $A^*$ (with $\lambda=0$).
	\item  The operator $\Gamma$ regularizes $\gamma_1$ in the following sense: $\Gamma u=0$ for $u\in\ker\At^*$, therefore
	$\Gamma u $ is determined only by the regular part of $u$ lying in $D(A_D)$. In fact we have:
\end{enumerate}
\end{remark}

\begin{theorem}[Grubb 1968]
Equip $D(\At^*)$ and $D(A^*)$ with the graph norm. Then $\Gamma:D(\At^*)\to H^{1/2}(\partial\Omega)$ is continuous and surjective. The same is true for 
$\Gammat:D(A^*)\to H^{1/2}(\partial\Omega)$. Moreover, for all $u\in D(\At^*)$, $v\in D(A^*)$ we have
\begin{eqnarray}\label{GreenGrubb}
(\At^* u, v)_{L^2(\Omega)} - (u,A^*v)_{L^2(\Omega)} = (\Gamma u, \gammat_2 v)_{\frac{1}{2},-\frac{1}{2}} - (\gamma_2 u, \Gammat v)_{-\frac{1}{2},\frac{1}{2}},
\end{eqnarray}
where $(\cdot,\cdot)_{\alpha,-\alpha}$ denotes the duality pairing between $H^\alpha(\partial\Omega)$ and $H^{-\alpha}(\partial\Omega)$.
\end{theorem}

\begin{proof}
See \cite[Theorem III.1.2]{Grubb}.
\end{proof}

To obtain an abstract Green formula of the form (\ref{Green}), we now need to rewrite the duality pairings as scalar products in $L^2(\partial\Omega)$. Since $L^2(\partial\Omega)$ and $H^{1/2}(\partial\Omega)$ are both infinite dimensional Hilbert spaces, there exists a unitary isomorphism $J:H^{1/2}(\partial\Omega)\to L^2(\partial\Omega)$. Then $(J^*)^{-1}:H^{-1/2}(\partial\Omega)\to L^2(\partial\Omega)$ is also a unitary isomorphism and 
\[(f, g)_{\frac{1}{2},-\frac{1}{2}}=(Jf,(J^*)^{-1}g)_{L^2(\partial\Omega)}.\]

\begin{theorem}
For $u\in D(\At^*)$ let \[\Gamma_1 u:=J\Gamma u,\quad \Gamma_2 u:=(J^*)^{-1}\gamma_2 u\]
and for $v\in D(A^*)$ let \[\Gammat_1 v:=J\Gammat v,\quad \Gammat_2 v:=(J^*)^{-1}\gammat_2 v.\]
Then
\[(\At^* u, v)_{L^2(\Omega)} - (u,A^*v)_{L^2(\Omega)} = (\Gamma_1 u, \Gammat_2 v)_{L^2(\partial\Omega)} - (\Gamma_2 u, \Gammat_1v)_{L^2(\partial\Omega)}.\]
Moreover,
\begin{enumerate}
	\item $\Gamma_i:D(\At^*)\to L^2(\partial\Omega)$ and $\Gammat_i:D(A^*)\to L^2(\partial\Omega)$ are surjective for $i=1,2$
	\item $\Gamma_i:D(\At^*)\to L^2(\partial\Omega)$ and $\Gammat_i:D(A^*)\to L^2(\partial\Omega)$ are continuous with respect to the graph norm for $i=1,2$
	\item given $(f,g)\in (L^2(\partial\Omega))^2$, there exist $u\in D(\At^*)$ such that $\Gamma_1 u=f$ and $\Gamma_2 u=g$ and $v\in D(A^*)$ such that $\Gammat_1 v=f$ and $\Gammat_2 v=g$ (inverse trace theorem).
\end{enumerate}
\end{theorem}

\begin{proof}
The Green identity follows from the previous theorem and the definition of $J$. 

Properties $(1)$ and $(2)$ are consequences of $\Gamma$ and $\Gammat$ being continuous and surjective onto $H^{1/2}(\partial\Omega)$ and $\gamma_2$ and $\gammat_2$ being continuous and surjective onto $H^{-1/2}(\partial\Omega)$ (\cf \cite[Proposition III.1.1]{Grubb}).

The inverse trace property $(3)$ follows from the corresponding property for $\Gamma$ and $\gamma_2$ and $\Gammat$ and $\gammat_2$, respectively  (\cf \cite[Proposition III.1.2]{Grubb}).
\end{proof}

\begin{remark}
\begin{itemize}
	\item All conditions we required in the previous sections on the boundary operators are satisfied by $\Gamma_1$, $\Gamma_2$, $\Gammat_1$ and $\Gammat_2$. So all the results on the corresponding $M$-functions hold.
	\item Note that $\At^*\big|_{\ker\Gamma_2}$ is the operator with Dirichlet boundary conditions - the Friedrichs extension of $A$, while $\At^*\big|_{\ker\Gamma_1}$ is the Kre\u{\i}n extension of $A$.
	\item By exchanging the roles of $\Gamma_1$ and $\Gamma_2$ it is possible to express the Neumann boundary condition in the form $\Gamma_1-B\Gamma_2$ for bounded $B$.
	\item An abstract form of this procedure for  regularizing the boundary operators has been introduced by Ryzhov \cite{Ryz07}.
\end{itemize}
\end{remark}

\end{document}